\documentclass[11pt, a4paper]{article}
\usepackage[dvips]{graphicx}
\usepackage[T1]{fontenc}
\usepackage[latin1]{inputenc}
\usepackage[english]{babel}
\usepackage{times}
\usepackage{amsmath}
\usepackage{amssymb}
\usepackage{amscd}
\usepackage{latexsym}
\usepackage{enumerate}
\usepackage{tabularx}
\usepackage{here}
\usepackage{color}
\usepackage{vmargin}
\usepackage{epsfig}
\usepackage{mathrsfs}

\usepackage{color}

\setmarginsrb{2cm}{2cm}{2cm}{2cm}{0cm}{2cm}{0cm}{1cm}

\begin{document}

\date{}

\title{DP-coloring for planar graphs of diameter two}

\author {\large Jingran Qi$^1$, Danjun Huang$^1$\footnote{This research was supported by Zhejiang Provincial Natural Science Foundation of China under Grant No. LY18A010014.}, Weifan Wang$^1$, Stephen Finbow$^2$ \\
{\small 1. Department of Mathematics,  Zhejiang Normal University, Jinhua, China}\\
{\small 2. Department of Mathematics and Statistics, St. Francis Xavier University, Antigonish, Canada}}

\maketitle

\newcommand{\qed}{\hfill $\Box$ }
\newcommand{\ch}{\omega}
\newcommand{\nch}{\omega^*}
\newtheorem{corollary}{Corollary}
\newtheorem{theorem}{Theorem}
\newtheorem{definition}{Definition}
\newtheorem{conjecture}{Conjecture}
\newtheorem{observation}{Observation}
\newtheorem{claim}{Claim}
\newtheorem{lemma}{Lemma}
\newtheorem{fact}{Fact}
\newtheorem{question}{Question}
\newtheorem{remark}{Remark}
\newcommand{\proof}{\noindent{\bf Proof.}\ \ }

\parindent=0.5cm

\begin{abstract}
\baselineskip 16pt
DP-coloring (also known as correspondence coloring) is a generalization of list coloring introduced by Dvo\u{r}\'{a}k and Postle (2017).  Recently, Huang et
al. [https://doi.org/10.1016/j.amc.2019.124562] showed that planar graphs with diameter at most two are $4$-choosable. In this paper, we will prove that planar graphs with diameter at most two are DP-$4$-colorable, which is an extension of the above result.

\medskip


\end{abstract}


\section{Introduction}

All graphs considered in this paper are finite, simple and undirected. Given a graph $G$ with vertex set $V(G)$
and edge set $E(G)$, a {\em proper $k$-coloring} of the graph $G$ is a mapping $f : V(G) \rightarrow \{1,2,\ldots,k\}$ such that $f(u) \neq
f(v)$ for any $uv \in E(G)$. The minimum integer $k$ such that $G$ admits a proper $k$-coloring is called the {\em
chromatic number} of $G$, denoted by $\chi(G)$.

A well-known generalization of proper $k$-coloring is the concept of list coloring, introduced by Vizing~\cite{Vizing-76}
and independently by Erd\H{o}s, Rubin, and Taylor~\cite{erd}.
We say that $L$ is a {\em $k$-list assignment} for the graph $G$ if it assigns a list $L(v)$ of possible colors to each
vertex $v$ of $G$ with $|L(v)|\ge k$. If $G$ has a proper coloring $f$ such that $f(v)\in L(v)$ for each vertex $v$, then we say that
$G$ is {\em $L$-colorable}. The graph $G$ is {\em $k$-choosable} if it is $L$-colorable for every $k$-list assignment $L$.
The {\em list chromatic number} of $G$, denoted by
$\chi_{l}(G)$, is the smallest positive integer $k$ such that $G$ is $k$-choosable.
By  the  definition, it holds trivially that $\chi_l(G)\ge \chi(G)$ for any graph $G$.

Thomassen \cite{T-1994} proved that every planar graph is 5-choosable, whereas Voigt \cite{Voigt-93}
presented an example of a planar graph which is not 4-choosable.
An interesting problem in graph coloring is to find sufficient
conditions for planar graphs to be 4-choosable.
It is easy to observe that planar graphs without 3-cycles are 4-choosable since they are 3-degenerate.
  Wang and Lih \cite{Wang-02} strengthened this result by showing that planar graphs without intersecting 3-cycles are 4-choosable, and
  conjectured that every planar graph without adjacent triangles is 4-choosable.
Moreover, several research groups proved that every planar graph without $k$-cycles  is 4-choosable  for $k=4$ in \cite{Lam-99}, for $k=5$ in \cite{Wang-01}, for $k=6$ in \cite{lam}, and for  $k=7$ in \cite{Farzad-09}.

Recently, Dvo\u{r}\'{a}k and Postle ~\cite{DP-2018} introduced
DP-coloring (under the name {\em correspondence coloring}) as a
generalization of list coloring.

\begin{definition} Let $G$ be a graph and $L$ be a list assignment of $V(G)$. For each edge $uv$ in $G$, let $M_{L,uv}$ be a
matching between the sets $\{u\} \times L(u)$ and $\{v\} \times L(v)$.
Let $\mathcal{M}_{L}=\{M_{L,uv} : uv \in E(G)\}$, which is called a {\em matching assignment} cover $L$.
Then a graph $H_L$ is said to be the $\mathcal{M}_{L}$-cover of $G$  if it satisfies all the following conditions:

{\rm (i)} The vertex set of $H_L$ is $\bigcup_{u \in V(G)}(\{u\} \times L(u))=\{(u,c) : u \in V(G), c \in L(u)\}$;

{\rm (ii)} For all $u \in V(G)$, the set $\{u\} \times L(u)$  induces a clique in $H_L$;

{\rm (iii)} If $uv \in E(G)$, then the edges between $\{u\} \times L(u)$ and $\{v\} \times L(v)$
            are those of $M_{L,uv}$.

{\rm (iv)} If $uv\notin E(G)$, then there are no edges between $\{u\} \times L(u)$ and $\{v\} \times L(v)$.
\end{definition}

\begin{definition}
An {\em $\mathcal{M}_{L}$-coloring} of $G$ is an independent set $I$ in the $\mathcal{M}_{L}$-cover with $|I|=|V(G)|$.
The {\em DP-chromatic number}, denoted by $\chi_{DP}(G)$, is the minimum integer $k$ such that $G$
admits an $\mathcal{M}_{L}$-coloring for each $k$-list assignment $L$ and each matching assignment $\mathcal{M}_{L}$ cover $L$. We say
that a graph $G$ is {\em DP-$k$-colorable} if $\chi_{DP}(G) \leq k$.
\end{definition}

As in list coloring, we refer to the elements of $L(v)$ as colors, and
call the element $i\in L(v)$ chosen in the independent set of an $\mathcal{M}_{L}$-coloring
as the color of $v$. And a vertex $v$ is colored with $i$ is also
stated that $i$ is the color of $v$.
Note that when $G$ is a simple graph with the list assignment $L$ and
$$M_{L,uv}=\{(u,c)(v,c):c \in L(u) \cap L(v)\}$$
for any edge $uv\in E(G)$, then $G$ admits an $L$-coloring if and only if $G$ admits an $\mathcal{M}_{L}$-coloring. This implies that
$\chi_{l}(G) \leq \chi_{DP}(G)$. We should note that DP-coloring and list coloring
can be quite different. For example, $\chi_l(C)=2<3=\chi_{DP}(C)$ for every even cycle $C$.

Dvo\u{r}\'{a}k and Postle~\cite{DP-2018} used Thomassen's proofs~\cite{T-1994} for choosability to
show that $\chi_{DP}(G) \leq 5$ for every planar graph $G$.
Recently, Kim and Ozeki ~\cite{KO-2018}
showed that planar graphs without $k$-cycles are DP-4-colorable for each $k=3,4,5,6$. Kim
and Yu ~\cite{KY-2019} extended the above results on $3$-cycles and $4$-cycles by showing that planar graphs without
triangles adjacent to $4$-cycles are DP-4-colorable.
Sittitrai and Nakprasit ~\cite{SN-2019} proved that every planar graph without pairwise adjacent $3$-,
$4$-, and $5$-cycle is DP-4-colorable. Some more sufficient conditions for a planar graph to be DP-4-colorable can
be found in \cite{CL-2019}, \cite{MO-2019}, \cite{LL-2019}.
In this paper, we show that every planar graph with diameter two is DP-$4$-colorable.


\section{Preliminaries}

We first introduce some notions used in this paper.
Suppose that $G$ is a planar graph embedded on the plane so that edges meet only at the vertices of the
graph. Let $F(G)$ be the set of faces. For a vertex $v\in V(G)$,  let $N(v)$ denote the set of neighbors of $v$,
and $d(v)=|N(v)|$.
We use $\delta(G)$ and $\Delta(G)$ to denote the minimum degree and the
maximum degree of $G$, respectively.
The {\em distance} between two vertices $u$ and $v$,
denoted by $d(u,v)$, is the length of a shortest path connecting them in $G$. The {\em diameter} of a graph $G$, denoted
diam$(G)$, is the maximum value of distances between any two vertices of $G$.

A planar graph is called {\em maximal} if $G+e$ is not planar for any edge $e\not\in E(G)$.
Clearly, if a maximal planar graph is embedded in the plane, then
each of its faces is a triangle.
For brevity,  a maximal planar graph of diameter two is called an
{\em MP$2$-graph}.
In 1989, Seyffarth ~\cite{S-1989} showed the following
results:

\begin{lemma}\label{maximal plane graph}
 Let $G$ be an MP$2$-graph. Then the following assertions
hold:

{\rm (1)} $3\leq \delta(G) \leq 4$.

{\rm (2)} If $\Delta(G) \geq 8$, then $|V(G)| \leq 1.5\Delta(G)+1$, and this bound is best possible.
\end{lemma}

Very recently, Huang et al.~\cite{HL-2019} characterized the structures of
every maximal planar graph $G$ of diameter two with $\delta(G)=4$.

\begin{lemma}\label{mp2-graph}
 Let $G$ be an MP$2$-graph with $\delta(G)=4$. Then $G\in \{ H_{1}, H_{2}, G_{1}, G_{2}, \ldots, G_{13}\}$, as shown
 in Figure.1.
\end{lemma}

\begin{figure}
  \centering
  \includegraphics[width=0.90\textwidth]{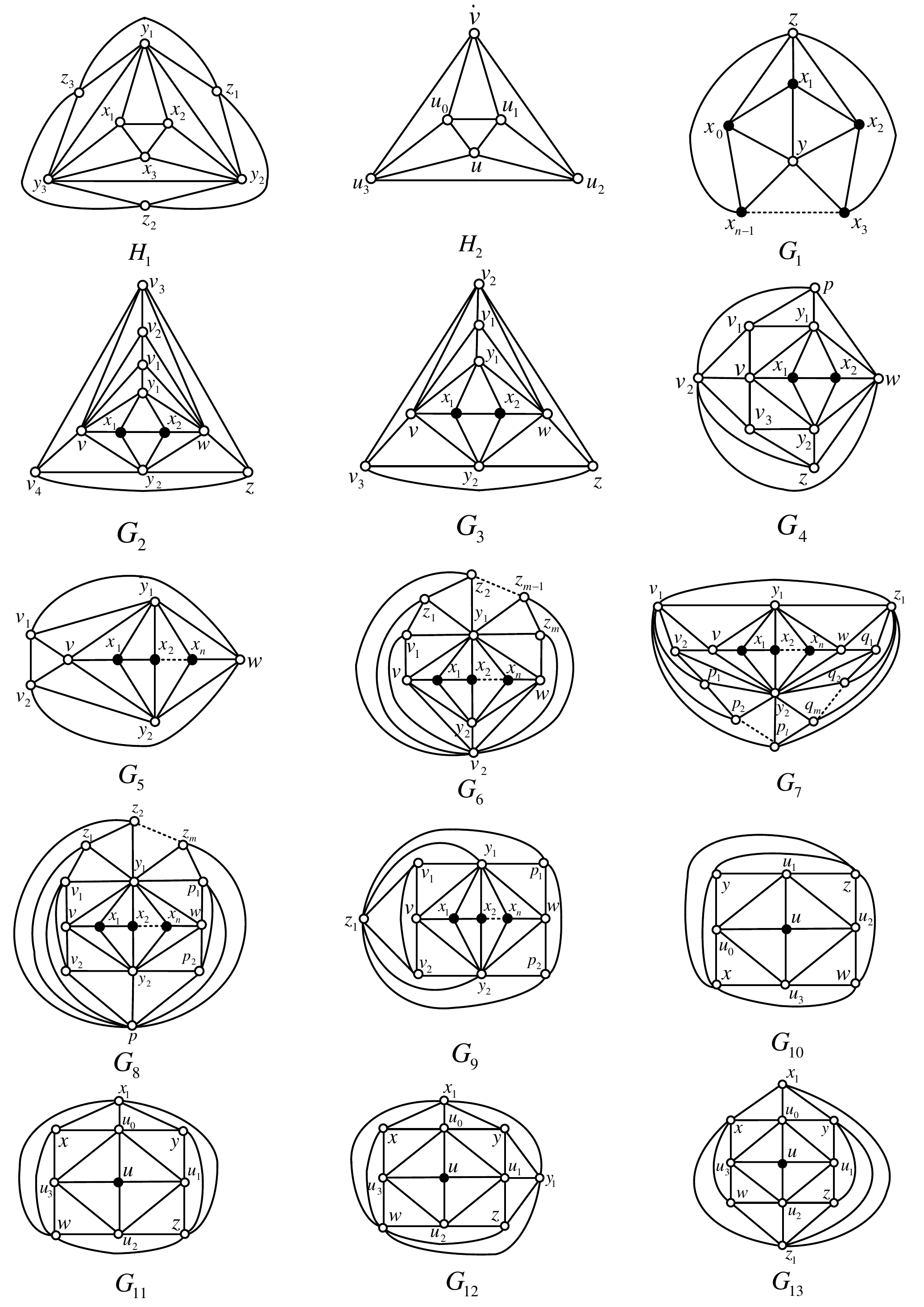}
  \caption{The structures of $G$ in Lemma~\ref{mp2-graph}. We redraw
  graphs $G_2$ and $G_3$ for the convenience of the proofs of Theorem~\ref{main2}. }\label{fig.1}
\end{figure}


\section{Main result and its proof}

Let $G$ be a plane graph, $L$ be a 4-list assignment of $G$,
and $\mathcal{M}_{L}$ be a matching assignment over $L$. An edge $uv\in E(G)$  is {\em straight} if
 every $(u,c_{1})(v,c_{2})\in M_{L,uv}$ satisfies $c_{1}=c_{2}$. The following lemma is from  (~\cite{DP-2018},
 Lemma 7) immediately.

 \begin{lemma}\label{sraight}
 Let $G$ be a graph with matching assignment $\mathcal{M}_{L}$. Let $H$ be a subgraph of $G$ which is a tree. Then we may
 rename $L(u)$ for $u\in  V(H)$ to obtain a matching assignment ${\mathcal{M}'_{L}}$ for $G$ such that all edges of $H$ are
 straight in ${\mathcal{M}'_{L}}$.
\end{lemma}

\begin{theorem}\label{main2}
Each MP$2$-graph $G$ is {\rm DP}-$4$-colorable.
\end{theorem}

\proof We will use induction on the vertex number of $G$ to prove Theorem~\ref{main2}. It is easy to
see that $G$ has at least five vertices. If $|V(G)|=5$, then $G$ is a graph
obtained from $K_{5}$ by removing one edge $e=uv$. Let $V(G)\setminus\{u,v\}=\{w_1,w_2,w_3\}$. By
Lemma~\ref{sraight}, we can rename each of $L(u),L(w_1),L(v)$ to
make $uw_1$ and $w_1v$ straight.  We can color $u$ and $v$ with the same color
$\alpha_{1}$, and then color $w_2$, $w_3$, $w_1$, successively.
So $G=K_5-e$ is DP-4-colorable.
Suppose that $G$ is an MP2-graph with $|V(G)|\geq 6$. By Lemma~\ref{maximal plane graph}, $3\leq \delta(G) \leq 4$. If
$\delta(G)=3$, say $d(v)=3$,  then we set $G^*=G-v$.  It is easy to check that $G^*$ is also an MP2-graph with
$|V(G^*)|=|V(G)|-1$.  By induction hypothesis, $G^*$ is  DP-$4$-colorable.
So there exists an independent set $I^*$ in $\mathcal{M}_{L}$-cover of $G^*$ with
$|I^*|=|V(G^*)|=|V(G)|-1$. Let $\phi(x)$ be the color of $x\in V(G^*)$ according to $I^*$. Obviously, we can extend
$I^*$ to the whole graph $G$ by coloring $v$ with a
color $\alpha$ such that for each $u\in N(v)$, $(u,\phi(u))(v,\alpha)\not\in M_{L,uv}$.
Hence, $G$ is DP-4-colorable.
So assume that $\delta(G)=4$. By Lemma~\ref{mp2-graph},  it suffices to testify that
$G$ is DP-$4$-colorable for each $G\in \{ H_{1}, H_{2}, G_{1}, G_{2}, \ldots, G_{13}\}$.

Let $L$ be a 4-list assignment of $G$, and
$\mathcal{M}_{L}$ be a matching assignment over $L$.
Without loss of generality, we may assume that $|L(v)|=4$ for each $v\in V(G)$, and
$M_{L,uv}$ is a perfect matching between the
sets $\{u\}\times L(u)$ and $\{v\}\times L(v)$ for each edge $uv\in E(G)$.
Let $H_L$ be the $\mathcal{M}_{L}$-cover of $G$.
We need to show that there exists an independent set $I$ in $H_L$ with $|I|=|V(G)|$.

Assume that $I'$ is an independent set in $H_L$ with
$|I'|<|V(G)|$.  For each $v\in V(G)$ with $(\{v\}\times L(v))\cap I'=\emptyset$,  define
$$L(v,I')=L(v)\setminus \bigcup_{uv\in E(G)} \{c'\in L(v) : (u,c)(v,c')\in M_{L,uv},\ \mbox{and}\ (u,c)\in I'\}.$$
Suppose that $C=v_{1}v_{2}v_{3}v_{1}$  is a $3$-cycle such that
$M_{L,v_{i}v_{i+1}}=\{(v_{i},c_{i,j})(v_{i+1},c_{i+1,j}) : j=1,2,3,4\}$ for $i=1,2$. We say that $C$ has {\em property $\mathcal{P}$} if
$M_{L,v_{1}v_{3}}=\{(v_{1},c_{1,j})(v_{3},c_{3,j}) : j=1,2,3,4\}$. That is, if
$C$ has property $\mathcal{P}$, then we can rename the colors in $L(v_{1})$, $L(v_{2})$ and $L(v_{3})$ to make $v_{1}v_{2}$,
$v_{2}v_{3}$ and $v_{3}v_{1}$ straight. If $C$ does not have property $\mathcal{P}$,
say $(v_{1},c_{1,1})(v_{3},c_{3,2})\in M_{L,v_{1}v_{3}}$, then we can color
$v_{1}$ with $c_{1,1}$ and color $v_{2}$ with
$c_{2,2}$ such that $|L(v_{3},I')|\geq 3$, where $I'=\{(v_1,c_{1,1}),(v_2,c_{2,2})\}$.

It is worth noting that if there exists a path $xyz$ with $xz\notin E(G)$ such that
$|L(x,I')|+|L(z,I')|>|L(y,I')|$ for some independent set $I'$ in $H_L$, then
we can color $x$ with $\alpha_1\in L(x,I')$ and color $z$ with $\alpha_2\in L(z,I')$ such that
$|L(y,I'')|\ge|L(y,I')|-1$, where $I''=I'\cup\{(x,\alpha_1),(z,\alpha_2)\}$.

\medskip
\noindent{\bf Case 1} $G=H_{1}$.
\medskip

Suppose that $C_{1}=y_{1}y_{2}y_{3}y_{1}$ does not have property $\mathcal{P}$.
Then there exist $\alpha_{1} \in L(y_{3})$ and $\alpha_{2} \in L(y_{2})$ such that
$|L(y_{1},I')|\geq 3$ when we color $y_{3}$ with
$\alpha_{1}$, color $y_{2}$ with $\alpha_{2}$, where $I'=\{(y_3,\alpha_1),(y_2,\alpha_2)\}$. Hence, $|L(z_{2},I')|\geq 2$
and $|L(z_{1},I')|\geq 3$. So we can color $z_{1}$ with $\alpha_{3} \in  L(z_{1},I')$
such that $|L(z_{2},I'')|\geq 2$, where $I''=I'\cup\{(z_1,\alpha_3)\}$. Similarly, we can color $x_{2}$ with $\alpha_{4}\in  L(x_{2},I'')$
such that $|L(x_{3},I''')|\geq 2$, where $I'''=I''\cup\{(x_2,\alpha_4)\}$.
Now we can color $y_{1}$, $x_{1}$, $x_{3}$, $z_{3}$, $z_{2}$, successively.

So suppose that $C_{1}=y_{1}y_{2}y_{3}y_{1}$ has the property $\mathcal{P}$.
That is, the edges $y_1y_2$, $y_2y_3$ and $y_3y_1$ are straight.

\begin{itemize}
\item Suppose that $C_{2}=x_{3}y_{2}y_{3}x_{3}$ does not have property $\mathcal{P}$.
Then there exist  $\alpha_{1} \in L(y_{3})$ and $\alpha_{2} \in L(y_{2})$ such that
$|L(x_{3},I')|\geq 3$ when we color $y_{3}$ with $\alpha_{1}$
and color $y_{2}$ with $\alpha_{2}$, where $I'=\{(y_3,\alpha_1),(y_2,\alpha_2)\}$.
Hence, $|L(y_{1},I')|\geq 2$, $|L(z_{1},I')|\geq 3$,
and $|L(z_{2},I')|\geq 2$. So we can color
$y_{1}$ with $\alpha_{3}\in L(y_{1},I')$ and color $z_{2}$ with $\alpha_{4}\in L(z_{2},I')$
such that $|L(z_{1},I'')|\geq 2$, where $I''=I'\cup\{(y_1,\alpha_3),(z_2,\alpha_4)\}$.
Now we can color $z_{3}$, $z_{1}$, $x_{1}$, $x_{2}$, $x_{3}$, successively.

\item Suppose that $C_{2}=x_{3}y_{2}y_{3}x_{3}$ has property $\mathcal{P}$.
That is, the edges $x_3y_2$, $y_2y_3$, and $y_3x_3$ are straight.
By symmetry, we may assume that $C_{3}=z_{2}y_{2}y_{3}z_{2}$,
$C_{4}=x_{2}y_{1}y_{2}x_{2}$, $C_{5}=z_{1}y_{1}y_{2}z_{1}$,
$C_{6}=x_{1}y_{1}y_{3}x_{1}$, and $C_{7}=z_{3}y_{1}y_{3}z_{3}$ have property $\mathcal{P}$.
We color $y_{2}, x_{1}$ and $z_{3}$ with the same color $\alpha_{1}$, and set $I'=\{(y_2,\alpha_1),(x_1,\alpha_1),(z_3,\alpha_1)\}$.
Hence, $|L(z_{1},I')|\geq 2$, $|L(y_{1},I')|\geq 3$, and $|L(x_{2},I')|\geq 2$.
So we can color $z_{1}$ with $\alpha_{2}\in L(z_{1},I')$ and color $x_{2}$ with
$\alpha_{3}\in  L(x_{2},I')$ such that $|L(y_{1},I'')|\geq 2$,
where $I''=I'\cup \{(z_1,\alpha_2),(x_2,\alpha_3)\}$.
Note that $|L(y_3,I'')|\geq 3$. Now we can color
$x_{3}$, $z_{2}$, $y_{3}$, $y_{1}$, successively.
\end{itemize}

\medskip
\noindent{\bf Case 2} $G=H_{2}$.
\medskip

Since $H_{2}$ is a subgraph of $H_{1}$, $H_{2}$ is DP-4-colorable.

\medskip
\noindent{\bf Case 3} $G=G_{1}$.
\medskip

By Lemma~\ref{sraight}, we can rename the colors in $L(y)$, $L(x_{1})$, $L(z)$ to
make $yx_{1}$ and $x_{1}z$ straight.
Now we can color $y$ and $z$ with the same color $\alpha_{1}$, and
color $x_{2}, x_{3}, \ldots, x_{n-1}, x_{0}, x_{1}$, successively.

\medskip
\noindent{\bf Case 4} $G=G_{2}$.
\medskip

Suppose that $C_{1}=vv_{2}v_{3}v$ does not have property $\mathcal{P}$.
Then we can color $v$ with $\alpha_{1}\in L(v)$ and color $v_{3}$ with $\alpha_{2} \in L(v_{3})$ such that
$|L(v_{2},I')|\geq 3$, where $I'=\{(v,\alpha_1),(v_3,\alpha_2)\}$.
Hence, $|L(v_{4},I')|\geq 2$ and $|L(y_{2},I')|\geq 3$.
So we can color $y_{2}$ with $\alpha_{3} \in  L(y_{2},I')$ such that $|L(v_{4},I'')|\geq 2$, where $I''=I'\cup\{(y_2,\alpha_3)\}$.
Similarly, we can color $y_{1}$
with $\alpha_{4}\in L(y_{1},I'')$ such that $|L(x_{1},I''')|\geq 2$, where $I'''=I''\cup\{(y_1,\alpha_4)\}$.
Now we can color $w$, $z$, $v_{4}$, $x_{2}$, $x_{1}$, $v_{1}$, $v_{2}$,  successively.

So suppose that $C_{1}=vv_{2}v_{3}v$ has property $\mathcal{P}$.
That is, the edges $vv_2$, $v_2v_3$ and $v_3v$ are straight.
By symmetry, $C_{2}=wv_{2}v_{3}w$ has property $\mathcal{P}$, too.
We can color $v$ and $w$ with the same color $\alpha_{1}$, and set
$I'=\{(v,\alpha_1),(w,\alpha_1)\}$. Hence, $|L(y_{1},I')|\geq 2$, $|L(x_{1},I')|\geq 3$, and $|L(y_{2},I')|\geq 2$.
So we can color $y_{1}$ with
$\alpha_{2}\in  L(y_{1},I')$ and color $y_{2}$ with $\alpha_{3}\in  L(y_{2},I')$,  such that
$|L(x_{1},I'')|\geq 2$, where $I''=I'\cup\{(y_1,\alpha_2),(y_2,\alpha_3)\}$.
Note that $|L(v_i,I'')|\geq 3$ for $i=2,3$.
Now we can  color $x_{2}$, $x_{1}$, $v_{1}$, $z$, $v_{4}$, $v_{3}$, $v_{2}$, successively.

\medskip
\noindent{\bf Case 5} $G=G_{3}$.
\medskip

Suppose that $C_{1}=vv_{1}v_{2}v$ does not have property $\mathcal{P}$.
Then we can color $v$ with $\alpha_{1}\in L(v)$ and color $v_{2}$ with $\alpha_{2}\in L(v_{2})$
such that $|L(v_{1},I')|\geq 3$, where $I'=\{(v,\alpha_1),(v_2,\alpha_2)\}$.
Hence, $|L(v_{3},I')|\geq 2$ and $|L(y_{2},I')|\geq 3$.
So we can color $y_{2}$ with $\alpha_{3}\in  L(y_{2},I')$ such that $|L(v_{3},I'')|\geq 2$, where $I''=I'\cup\{(y_2,\alpha_3)\}$.
Similarly, we can color $y_{1}$ with $\alpha_{4}\in  L(y_{1},I'')$ such that
$|L(x_{1},I''')|\geq 2$, where $I'''=I''\cup\{(y_1,\alpha_4\}$.
Now we can color $w$, $v_{1}$, $x_{2}$, $x_{1}$, $z$, $v_{3}$, successively.

So suppose that $C_{1}=vv_{1}v_{2}v$ has property $\mathcal{P}$.
That is, the edges $vv_1$, $v_1v_2$ and $v_2v$ are straight.
By symmetry, $C_{2}=wv_{1}v_{2}w$ has property $\mathcal{P}$, too.
We can color $v$ and $w$ with the same color $\alpha_{1}$, and set
$I'=\{(v,\alpha_1),(w,\alpha_1)\}$. Hence, $|L(y_{1},I')|\geq 2$, $|L(x_{1},I')|\geq 3$, and $|L(y_{2},I')|\geq 2$.
So we can color $y_{1}$ with
$\alpha_{2}\in  L(y_{1},I')$ and color $y_{2}$ with $\alpha_{3}\in  L(y_{2},I')$,  such that
$|L(x_{1},I'')|\geq 2$, where $I''=I'\cup\{(y_1,\alpha_2),(y_2,\alpha_3)\}$.
Note that $|L(v_1,I'')|\geq 2$ and $|L(v_2,I'')|\geq 3$.
Now we can  color $x_{2}$, $x_{1}$, $z$, $v_{3}$, $v_{2}$, $v_{1}$, successively.

\medskip
\noindent{\bf Case 6} $G=G_{4}$.
\medskip

Suppose that $C_{1}=vy_{1}x_{1}v$ does not have property $\mathcal{P}$.
Then we can color $v$ with $\alpha_{1}\in L(v)$ and color $y_{1}$ with $\alpha_{2}\in L(y_{1})$ such that
$|L(x_{1},I')|\geq 3$, where $I'=\{(v,\alpha_1),(y_1,\alpha_2)\}$.
Hence, $|L(v_{1},I')|\geq 2$ and $|L(v_{2},I')|\geq 3$.
So we can color $v_2$ with $\alpha_{3}\in  L(v_{2},I')$ such that $|L(v_{1},I'')|\geq 2$, where $I''=I'\cup\{(v_2,\alpha_3)\}$.
Similarly, we can color $y_2$ with $\alpha_{4}\in  L(y_{2},I'')$ such that
$|L(v_{3},I''')|\geq 2$, where $I'''=I''\cup\{(y_2,\alpha_4)\}$.
Now we can color $w$, $x_{2}$, $x_{1}$, $z$, $v_{3}$, $p$, $v_{1}$, successively.

So suppose that $C_{1}=vy_{1}x_{1}v$ has property $\mathcal{P}$.
That is, the edges $vy_1$, $y_1x_1$ and $x_1v$ are straight.
By symmetry, $C_{2}=vy_{2}x_{1}v$ has property $\mathcal{P}$, too.
By Lemma~\ref{sraight},
we can rename the colors in $L(y_{1})$, $L(v_{1})$
and $L(v_{2})$ to make $y_{1}v_{1}$ and $v_{1}v_{2}$ straight.
We can color $y_{1}$, $y_2$ and $v_{2}$ with the same color $\alpha_{1}$,
and set $I'=\{(y_1,\alpha_1),(y_2,\alpha_1),(v_2,\alpha_1)\}$.
Note that $|L(v_{1},I')|\geq 3$, $|L(v,I')|\geq 2$ and $|L(x_1,I')|\geq 3$.
Now we color $w$, $x_{2}$, $z$, $v_{3}$, $v$, $x_{1}$, $p$, $v_{1}$, successively.

\medskip
\noindent{\bf Case 7} $G=G_{5}$.
\medskip

By Lemma~\ref{sraight}, we can rename the colors in $L(y_{1})$, $L(x_{1})$ and $L(y_{2})$ to
make $y_{1}x_{1}$ and $x_{1}y_{2}$ straight.
Then we can color $y_{1}$ and $y_{2}$ with the same color $\alpha_{1}$, and set $I'=\{(y_1,\alpha_1),(y_2,\alpha_1)\}$.
Hence, $|L(v,I')|\geq 2$, $|L(v_{2},I')|\geq 3$ and $|L(w,I')|\geq 2$.
So we can color $v$ with $\alpha_{2}\in  L(v,I')$ and color $w$ with $\alpha_{3}\in  L(w,I')$  such that
$|L(v_{2},I'')|\geq 2$, where $I''=I'\cup\{(v,\alpha_2),(w,\alpha_3)\}$.
Note that $|L(x_1,I'')|\ge2$. Now we can color $x_{n},\ldots, x_{2}$, $x_{1}$, $v_{1}$, $v_{2}$, successively.

\medskip
\noindent{\bf Case 8} $G=G_{6}$.
\medskip

By Lemma~\ref{sraight}, we can rename the colors in $L(y_{1})$, $L(x_{1})$ and $L(y_{2})$ to
make $y_{1}x_{1}$ and $x_{1}y_{2}$ straight.
Then we color $y_{1}$ and $y_{2}$ with the same color $\alpha_{1}$, and set $I'=\{(y_1,\alpha_1),(y_2,\alpha_1)\}$.
Hence, $|L(w,I')|\geq 2$ and $|L(v_{2},I')|\geq 3$.
So we can color $v_{2}$ with $\alpha_{2}\in L(v_{2},I')$ such that $|L(w,I'')|\geq 2$,
where $I''=I'\cup\{(v_2,\alpha_2)\}$. Note that $|L(x_1,I'')|\ge3$.
Now we can color $v$, $v_{1}$, $z_{1}$,
$z_{2},\ldots, z_{m}$, $w$, $x_{n},\ldots, x_{2}$, $x_{1}$, successively.

\medskip
\noindent{\bf Case 9} $G=G_{7}$.
\medskip

Suppose that $C_{1}=p_{l}y_{2}p_{l-1}p_{l}$ does not have property $\mathcal{P}$.
Then we can color $p_{l}$ with $\alpha_{1}\in L(p_{l})$ and
color $y_{2}$ with $\alpha_{2}\in L(y_{2})$ such that $|L(p_{l-1},I')|\geq 3$, where $I'=\{(p_l,\alpha_1),(y_2,\alpha_2)\}$.
Hence, $|L(q_{m},I')|\geq 2$ and $|L(z_{1},I')|\geq 3$.
So we can color $z_{1}$ with $\alpha_{3}\in  L(z_{1},I')$ such that
$|L(q_{m},I'')|\geq 2$, where $I''=I'\cup\{z_1,\alpha_3)\}$.
Similarly, we can color $y_{1}$ with $\alpha_{4}\in L(y_{1},I'')$ such that
$|L(w,I''')|\geq 2$, where $I'''=I''\cup\{(y_1,\alpha_4)\}$.
Now we can  color $v_{1}$, $v$, $x_{1}$, $x_{2},\ldots, x_{n}$, $w$, $q_{1}$, $q_{2},\ldots, q_{m}$,
$v_{2}$, $p_{1}$, $p_{2},\ldots, p_{l-1}$, successively.

So suppose that $C_{1}=p_{l}y_{2}p_{l-1}p_{l}$ has property $\mathcal{P}$.
That is, $p_{l}p_{l-1}$, $p_{l}y_{2}$ and $p_{l-1}y_{2}$ are  straight.
If $C_{2}=v_{1}p_{l}p_{l-1}v_{1}$ has property $\mathcal{P}$,
say $v_1p_l$, $p_lp_{l-1}$ and $p_{l-1}v_1$ are straight,
then we can color $v_{1}$ and $y_{2}$ with the same color $\alpha_{1}$ such that
$|L(p_{l},I')|\geq 3$ and $|L(p_{l-1},I')|\geq 3$, where $I'=\{(v_1,\alpha_1),(y_2,\alpha_1)\}$.
Hence, $|L(v,I')|\geq 2$ and $|L(y_{1},I')|\geq 3$.
So we can color then $y_{1}$ with $\alpha_{2}\in  L(y_{1},I')$ such that
$|L(v,I'')|\geq 2$, where $I''=I'\cup\{(y_1,\alpha_2)\}$.
Now we can color $z_{1}$,
$w$, $x_{n},\ldots, x_{1}$, $v$, $v_{2}$, $p_{1}$, $p_{2},\ldots, p_{l-2}$, $q_{1}$,
$q_{2},\ldots, q_{m}$, $p_{l}$, $p_{l-1}$, successively. Hence, we may assume that
$C_{2}=v_{1}p_{l}p_{l-1}v_{1}$ does not have property $\mathcal{P}$.

From now on, we  rename the colors in $L(y_{1})$, $L(v)$, $L(v_{1})$, $L(y_{2})$
and $L(w)$ to make $y_{1}v$, $vv_{1}$, $vy_{2}$ and $y_{2}w$ straight by Lemma~\ref{sraight}.
Recall that $p_{l}p_{l-1}$, $p_{l}y_{2}$ and $p_{l-1}y_{2}$ are  straight.
Since $C_{2}=v_{1}p_{l}p_{l-1}v_{1}$ does not have property $\mathcal{P}$,
we can color $v_{1}$ with $\alpha_{1}\in L(v_1)$ and $p_{l}$ with $\alpha_{2}\in L(p_l)$
such that $|L(p_{l-1},I')|\geq 3$, where $I'=\{(v_1,\alpha_1),(p_l,\alpha_2)\}$.

\begin{itemize}
\item Suppose that $\alpha_1=\alpha_2$. Then $(p_{l-1},\alpha_1)(p_l,\alpha_1)\in M_{L,p_{l-1}p_l}$
by $p_lp_{l-1}$ is straight. So $(v_{1},\alpha_{1})(p_{l-1},\alpha_{1})\in M_{L,v_{1}p_{l-1}}$ by $|L(p_{l-1},I')|\geq 3$.
Now we remove the color of $p_l$, and color $y_{2}$ with
$\alpha_{1}$. Set $I''=\{(v_1,\alpha_1),(y_2,\alpha_1)\}$.
It is trivial that $|L(p_{l-1},I'')|\geq 3$ and $|L(v,I'')|\geq 3$.
Hence, $|L(p_{l},I'')|\geq 2$ and
$|L(z_{1},I'')|\geq 3$. So we can color $z_1$ with $\alpha_{3}\in  L(z_{1},I'')$ such that
$|L(p_{l},I''')|\geq 2$, where $I'''=I''\cup\{(z_1,\alpha_3)\}$.
Now we color $y_{1}$, $w$, $x_{n},\ldots, x_{1}$, $v$, $v_{2}$, $p_{1}$, $p_{2},\ldots,
p_{l-2}$, $q_{1}$, $q_{2},\ldots, q_{m}$, $p_{l}$, $p_{l-1}$, successively.

\item Suppose that $\alpha_1\neq\alpha_2$. We
color $y_{2}$ with $\alpha_{1}$, and set $I''=I'\cup\{(y_2,\alpha_1)\}$.
It is trivial that $|L(v,I'')|\geq 3$ and $|L(p_{l-1},I'')|\ge2$.
Now we color $z_1$, $q_{m},\ldots, q_{1}$,
$w$, $y_1$, $x_{n},\ldots, x_{1}$, $v$, $v_{2}$, $p_{1}$, $p_{2},\ldots, p_{l-1}$, successively.
\end{itemize}

\medskip
\noindent{\bf Case 10} $G=G_{8}$.

\begin{claim}\label{py2v2}
If $C_{1}=py_{2}v_{2}p$ does not have property $\mathcal{P}$, then there is an independent set $I$ in $\mathcal{M}_L$-cover of $G_8$, say $H_L$, with $|I|=|V(G_8)|$.
\end{claim}

\proof Then we can color $p$ with $\alpha_{1}\in L(p)$
 and color $y_{2}$ with $\alpha_{2}\in L(y_2)$ such that $|L(v_{2},I')|\geq 3$, where $I'=\{(p,\alpha_1), (y_2,\alpha_2)\}$.
 Hence, $|L(p_{2},I')|\geq 2$ and $|L(w,I')|\geq 3$.
So we can color $w$ with $\alpha_{3}\in L(w,I')$ such that $|L(p_{2},I'')|\geq 2$, where $I''=I'\cup\{(w,\alpha_3)\}$.
Similarly, we can color $y_{1}$
 with $\alpha_{4}\in L(y_1,I'')$ such that $|L(x_{n},I''')|\geq 2$, where $I'''=I''\cup\{(y_1,\alpha_4)\}$.
 Now we can color $p_{1}$, $p_{2}$, $z_{m},\ldots, z_{1}$, $v_{1}$, $v$,
 $v_{2}$, $x_{1},\ldots, x_{n}$, successively.\qed

 \begin{claim}\label{y1xnwy2}
If there is $\alpha_{1} \in L(y_{1})$ and $\alpha_2\in L(y_2)$ such that $|L(x_i,I')|\ge3$ and $|L(w,I')|\ge3$ (or $|L(v,I')|\ge3$), where
$I'=\{(y_1,\alpha_1),(y_2,\alpha_2)\}$ and $i\in\{1,2,\ldots,n\}$, then there is an independent set $I$ in $H_L$ with $|I|=|V(G_8)|$.
 \end{claim}

 \proof By symmetry, we may assume that $|L(w,I')|\ge3$. We first color $y_{1}$ with $\alpha_1$ and color $y_{2}$ with $\alpha_{2}$.
 Hence, $|L(v,I')|\geq 2$ and $|L(v_{2},I')|\geq 3$.
 So we can color $v_{2}$ with $\alpha_{3}\in L(v_2,I')$ such that  $|L(v,I'')|\geq 2$, where
 $I''=I'\cup\{(v_2,\alpha_3)\}$. Now we color $p$, $v_{1}$, $v$,
 $x_{1},\ldots, x_{i-1}$, $z_{1},\ldots, z_{m}$, $p_{1}$, $p_{2}$, $w$, $x_{n},\ldots,x_i$, successively.
\qed

We can show the following claim analogously.

\begin{claim}\label{y1p1pzm}
If there is $\alpha_{1} \in L(y_{1})$ and $\alpha_2\in L(p)$ such that $|L(z_i,I')|\ge3$ and $|L(p_1,I')|\ge3$ (or $|L(v_1,I')|\ge3$), where
$I'=\{(y_1,\alpha_1),(p,\alpha_2)\}$ and $i\in\{1,2,\ldots,m\}$, then there is an independent set $I$ in $H_L$ with $|I|=|V(G_8)|$.
 \end{claim}

 \begin{claim}\label{y1vwy2}
If there is $\alpha_{1} \in L(y_{1})$ and $\alpha_2\in L(y_2)$ such that $|L(v,I')|\ge3$ and $|L(w,I')|\ge3$, where
$I'=\{(y_1,\alpha_1),(y_2,\alpha_2)\}$, then there is an independent set $I$ in $H_L$ with $|I|=|V(G_8)|$.
 \end{claim}

 \proof We first color $y_{1}$ and $\alpha_1$ and color $y_{2}$ with $\alpha_{2}$.
Hence, $|L(x_{n},I')|\geq 2$. So we can color $w$ with $\alpha_{3}$ such that $|L(x_{n},I'')|\geq 2$, where $I''=I'\cup\{(w,\alpha_3)\}$.
Similarly, we can color $p$ with $\alpha_{4}$ such that $|L(p_{1},I''')|\geq 2$, where $I'''=I''\cup\{(p,\alpha_4)\}$.
Now we can color $p_{2}$,
 $p_{1}$, $z_{m},\ldots, z_{1}$, $v_{1}$, $v_{2}$, $v$, $x_{1},\ldots, x_{n}$, successively.
\qed

Similarly, we can obtain the following claim.

\begin{claim}\label{y1v1p1p}
If there is $\alpha_{1} \in L(y_{1})$ and $\alpha_2\in L(p)$ such that $|L(v_1,I')|\ge3$ and $|L(p_1,I')|\ge3$, where
$I'=\{(y_1,\alpha_1),(p,\alpha_2)\}$, then there is an independent set $I$ in $H_L$ with $|I|=|V(G_8)|$.
 \end{claim}

\noindent{\bf Remark}. Suppose that $(y_1,\alpha_1)(p_1,\alpha_1)\in M_{L,y_1p_1}$,
$(y_1,\alpha_1)(z_m,\alpha_1)\in M_{L,y_1z_m}$ and $(y_1,\alpha_1)(v_1,\alpha_1)\in M_{L,y_1v_1}$.
Let $(p_1,\alpha_1)(p,\alpha)\in M_{L,p_1p}$, where $\alpha\in L(p)$.
If $(z_m,\alpha_1)(p,\alpha)\in M_{L,z_mp}$,
then $|L(p_1,\{(y_1,\alpha_1),(p,\alpha)\})|\ge3$ and $|L(z_m,\{(y_1,\alpha_1),(p,\alpha)\})|\ge3$.
So there exists the required independent set by Claim~\ref{y1p1pzm}. Hence,
$(z_m,\alpha_1)(p,\alpha)\notin M_{L,z_mp}$. Similarly, we can show
that $(p,\alpha)(v_{1},\alpha_{1}) \notin M_{L,pv_{1}}$ by Claim~\ref{y1v1p1p}.
\medskip

Now we are ready to show that there is an independent set $I$ in $H_L$ with $|I|=|V(G_8)|$.
 That is, $G_8$ is DP-4-colorable.

\medskip
\noindent{\bf Case 10.1} $C_{2}=vy_{2}v_{2}v$ does not have property $\mathcal{P}$.
\medskip

So we can color $y_{2}$ with
 $\alpha_{1}\in L(y_2)$ and color $v$ with $\alpha_{2}\in L(v)$ such that $|L(v_{2},I_1)|\geq 3$,
where $I_1=\{(y_2,\alpha_1),(v,\alpha_2)\}$.
If $C_1=py_2v_2p$ does not have property $\mathcal{P}$, then
there exists the required independent set $I$ by Claim~\ref{py2v2}. So suppose that
$C_1$ has property $\mathcal{P}$. By symmetry, $py_2p_2p$ has property
$\mathcal{P}$, too. By Lemma~\ref{sraight}, we can rename the colors in $L(x)$, where $x\in V(G_8)$, to make
$y_{1}v_{1}$, $y_{1}z_{m}$, $y_{1}p_{1}$, $y_{1}w$, $wy_{2}$, $y_{2}x_{n}$, $y_{2}v$, $y_2v_2$,
$y_2p_2$ and $y_2p$ straight. Then the edges $v_2p$ and $p_2p$ are also straight, since both $py_2v_2p$
and $py_2p_2p$ have property $\mathcal{P}$.

\medskip
\noindent{\bf Case 10.1.1} $(v,\alpha_{2})(y_{1},\alpha_{1})\not \in M_{L,vy_{1}}$.
\medskip

Then color $y_{1}$ with $\alpha_{1}$, and set $I_2=I_1\cup\{(y_1,\alpha_1)\}$.
Note that $y_1w$ and $wy_2$ are straight. So $|L(w,I_2)|\geq 3$.
If $|L(x_n,I_2)|\ge3$, then there exists the required independent set $I$ by Claim~\ref{y1xnwy2}.
So $|L(x_n,I_2)|=2$, say $(y_1,\alpha_1)(x_n,\beta_1)\in M_{L,y_1x_n}$, where $\beta_1\neq\alpha_1$.

First, suppose that $(w,\alpha_{1})(p_{2},\alpha_{1})\in M_{L,wp_{2}}$.
If $(p_{1},\alpha_{1})(p_{2},\alpha_{1})\in M_{L,p_{1}p_{2}}$, then
we remove the colors of $y_{2}$ and $v$, color $p_{2}$ with $\alpha_{1}$,
and color $y_2$ with $\beta_1$.
Set $I_3=\{(y_2,\beta_1),(y_1,\alpha_1),(p_2,\alpha_1)\}$.
Then $|L(p_{1},I_3)|\geq 3$, $|L(x_n,I_3)|\ge3$ and $|L(w,I_3)|\geq 2$ by the assumption.
Note that $|L(v_2,I_3)|\ge3$ and $|L(p,I_3)|\ge2$.
So we can color $v_{2}$ with $\alpha_3\in L(v_2,I_3)$ such that $|L(p,I_3\cup\{(v_2,\alpha_3)\})|\geq 2$.
Now we can color $v$, $v_{1}$, $p$, $z_{1},\ldots, z_{m}$, $p_{1}$, $w$, $x_{1},\ldots, x_{n}$, successively.
So we may suppose that $(p_{1},\alpha_{1})(p_{2},\alpha_{1})\notin M_{L,p_{1}p_{2}}$. Then
$(p_{1},\alpha_{1})(p_{2},\beta_2)\in M_{L,p_{1}p_{2}}$, where $\beta_2\neq \alpha_{1}$.
We color $p_{2}$ with $\beta_2$ and set $I'_1=I_2\cup\{(p_2,\beta_2)\}$.
Now $|L(p_{1},I'_1)|\geq 3$ and $|L(w,I'_1)|\ge2$. Suppose that $(p_1,\alpha_1)(p,\alpha)\in M_{L,p_1p}$, where $\alpha\in L(p)$.
If $\alpha\notin \{\alpha_1,\beta_2\}$, then we color $p$ with $\alpha$ and set
$I'_2=I'_1\cup\{(p,\alpha)\}$. Hence, $|L(p_1,I'_2)|\ge3$ and $|L(v_2,I'_2)|\ge2$.
Now we can color $v_{1}$, $v_{2}$, $z_{1},\ldots, z_{m}$, $x_{1},\ldots, x_{n}$, $w$, $p_1$, successively.
Otherwise, $\alpha\in\{\alpha_1,\beta_2\}$. Let $(z_{m},\alpha_{1})(p,\beta') \in M_{L,pz_{m}}$.
Then $\beta'\neq\alpha$ by the Remark. If $\beta' \notin
\{\alpha_{1},\beta_2\}$, then we color $p$ with $\beta'$. Hence,  $|L(z_{m},I'_1\cup\{(p,\beta')\})|\geq 3$. Now we can color
$v_{1}$, $v_{2}$, $z_{1},\ldots, z_{m-1}$, $x_{1},\ldots, x_{n}$, $w$, $p_{1}$,
$z_{m}$, successively. So suppose that $\beta'\in\{\alpha_{1},\beta_2\}\setminus\{\alpha\}$.
Recall that $(p,\alpha)(v_{1},\alpha_{1}) \notin M_{L,pv_{1}}$ by the Remark.
If $(p,\beta')(v_{1},\alpha_{1}) \in M_{L,pv_{1}}$,
then $|L(v_{1},\{(y_1,\alpha_1),(p,\beta')\})|\geq 3$ and $|L(z_{m},\{(y_1,\alpha_1),(p,\beta')\})|\geq 3$.
By Claim~\ref{y1p1pzm}, there exists the required independent set.
Hence, $(p,\beta')(v_{1},\alpha_{1})\notin M_{L,pv_{1}}$, say $(p,\gamma)(v_{1},\alpha_{1})\in M_{L,pv_{1}}$,
where $\gamma\notin\{\alpha_{1},\beta_2\}$. We color
 $p$ with $\gamma$. Then $|L(u,I'_1\cup\{(p,\gamma)\})|\geq 2$, where $u\in\{v_1,v_2,p_1,w\}$.
 Now we can color $x_{1},\ldots, x_{n}$, $w$, $p_{1}$, $z_{m},\ldots, z_{1}$, $v_{1}$, $v_{2}$,
 successively.

Next, suppose that $(w,\alpha_{1})(p_{2},\alpha_{3})\in M_{L,wp_{2}}$ with $\alpha_{3}\neq \alpha_{1}$.
We color $p_2$ with $\alpha_3$ and set $I''_1=I_2\cup\{(p_2,\alpha_3)\}$.
Now $|L(w,I''_1)|\ge3$.
Suppose that $(p_1,\alpha_1)(p,\alpha)\in M_{L,p_1p}$, where $\alpha\in L(p)$.
If $\alpha\notin \{\alpha_{1},\alpha_{3}\}$, then we color $p$ with $\alpha$.
Hence, $|L(p_{1},I''_1\cup\{(p,\alpha)\})|\geq 2$. Now we can color $v_{1}$, $v_{2}$, $z_{1},\ldots, z_{m}$,
$p_{1}$, $x_{1},\ldots, x_{n}$, $w$, successively. So suppose that $\alpha\in \{\alpha_{1},\alpha_{3}\}$.
Let $(z_{m},\alpha_{1})(p,\beta'') \in M_{L,pz_{m}}$. Then $\beta''\neq\alpha$ by the Remark.
If $\beta''\notin \{\alpha_1,\alpha_{3}\}$, then we color
$p$ with $\beta''$. Hence, $|L(z_{m},I''_1\cup\{(p,\beta'')\})|\geq 3$. Now we can color $v_{1}$, $v_{2}$, $z_{1},\ldots, z_{m-1}$, $p_{1}$,
 $z_{m}$, $x_{1},\ldots, x_{n}$, $w$, successively. So suppose that $\beta''\in\{\alpha_1,\alpha_{3}\}\setminus\{\alpha\}$.
Recall that $(p,\alpha)(v_{1},\alpha_{1}) \notin M_{L,pv_{1}}$ by the Remark.
If $(p,\beta'')(v_{1},\alpha_{1}) \in M_{L,pv_{1}}$,
then $|L(v_{1},\{(y_1,\alpha_1),(p,\beta'')\})|\geq 3$ and $|L(z_{m},\{(y_1,\alpha_1),(p,\beta'')\})|\geq 3$.
By Claim~\ref{y1p1pzm}, there exists the required independent set.
 Hence, $(p,\beta'')(v_{1},\alpha_{1}) \notin M_{L,pv_{1}}$, say $(p,\gamma')(v_{1},\alpha_{1}) \in M_{L,pv_{1}}$,
 where $\gamma' \notin \{\alpha_{1},\alpha_{3}\}$.
We color
 $p$ with $\gamma'$. Then $|L(v_1,I'_1\cup\{(p,\gamma)\})|\geq 2$ and $|L(v_2,I'_1\cup\{(p,\gamma)\})|\geq 2$.
Now we can color $p_{1}$, $z_{m},\ldots, z_{1}$, $v_{1}$, $v_{2}$,
 $x_{1},\ldots, x_{n}$, $w$, successively.

\medskip
\noindent{\bf Case 10.1.2} Suppose that $(v,\alpha_{2})(y_{1},\alpha_{1})\in M_{L,vy_{1}}$.
\medskip

Then $|L(w,\{(y_1,\alpha_1),(y_2,\alpha_1)\}|\ge3$. If $|L(x_n,\{(y_1,\alpha_1),(y_2,\alpha_1)\}|\ge3$,
then there exists the required independent set by Claim~\ref{y1xnwy2}. So $|L(x_n,\{(y_1,\alpha_1),(y_2,\alpha_1)\}|=2$.
Hence, we can color $y_{1}$ with ${\alpha_{3}}$ such that $|L(x_{n},I_2)|\geq 3$,
where ${\alpha_{3}} \neq \alpha_{1}$ and $I_2=I_1\cup\{(y_1,\alpha_3)\}$.
Therefore, $|L(w,I_2)|\geq 2$ and $|L(p_{2},I_2)|\geq 3$. Then
we can color $p_{2}$ with $\alpha_{4}\in L(p_2,I_2)$ such that $|L(w,I_3)|\geq 2$, where
$I_3=I_2\cup\{(p_2,\alpha_4)\}$. Suppose that $(p_1,\alpha_3)(p,\alpha)\in M_{L,p_1p}$, where $\alpha\in L(p)$.
If $\alpha\notin \{\alpha_1,\alpha_4\}$, then we color $p$ with $\alpha$ and set
$I_4=I_3\cup\{(p,\alpha)\}$. Hence, $|L(p_1,I_4)|\ge2$ and $|L(v_2,I_4)|\ge2$.
Now we can color $v_{1}$, $v_{2}$, $z_{1},\ldots, z_{m}$, $p_1$, $w$, $x_{1},\ldots, x_{n}$, successively.
Otherwise, $\alpha\in\{\alpha_1,\alpha_4\}$. Let $(z_{m},\alpha_{3})(p,\beta) \in M_{L,pz_{m}}$.
Then $\beta\neq\alpha$ by the Remark. If $\beta\notin
\{\alpha_{1},\alpha_4\}$, then we color $p$ with $\beta$. Hence,  $|L(z_{m},I_3\cup\{(p,\beta)\})|\geq 3$. Now we can color
$v_{1}$, $v_{2}$, $z_{1},\ldots, z_{m-1}$, $p_{1}$, $z_{m}$, $w$, $x_{1},\ldots, x_{n}$,
successively. So suppose that $\beta\in\{\alpha_{1},\alpha_4\}\setminus\{\alpha\}$.
Recall that $(p,\alpha)(v_{1},\alpha_{3}) \notin M_{L,pv_{1}}$ by the Remark.
If $(p,\beta)(v_{1},\alpha_{3}) \in M_{L,pv_{1}}$,
then $|L(v_{1},\{(y_1,\alpha_3),(p,\beta)\})|\geq 3$ and $|L(z_{m},\{(y_1,\alpha_3),(p,\beta)\})|\geq 3$.
By Claim~\ref{y1p1pzm}, there exists the required independent set.
Hence, $(p,\beta)(v_{1},\alpha_{3})\notin M_{L,pv_{1}}$, say $(p,\gamma)(v_{1},\alpha_{3})\in M_{L,pv_{1}}$,
where $\gamma\notin\{\alpha_{1},\alpha_4\}$. We color
$p$ with $\gamma$. Then $|L(u,I_3\cup\{(p,\gamma)\})|\geq 2$, where $u\in\{v_1,v_2,w\}$, and $|L(x_n,I_3\cup\{(p,\gamma)\})|\geq 3$.
Now we can color $p_{1}$, $z_{m},\ldots, z_{1}$, $v_{1}$, $v_{2}$, $w$, $x_{1},\ldots, x_{n}$,
 successively.

\medskip
\noindent{\bf Case 10.2} $C_{2}=vy_{2}v_{2}v$ has property $\mathcal{P}$.
\medskip

By symmetry, $C_{3}=wy_{2}p_{2}w$ has the property $\mathcal{P}$, too.
If $C_1=py_2v_2p$ does not have property $\mathcal{P}$, then
there exists the required independent set $I$ by Claim~\ref{py2v2}. So suppose that
$C_1$ has property $\mathcal{P}$. By symmetry, $py_2p_2p$ has property
$\mathcal{P}$, too. By Lemma~\ref{sraight}, we can rename the colors in $L(x)$, where $x\in V(G_8)$, to make
$y_{1}v_{1}$, $y_{1}z_{m}$, $y_{1}p_{1}$, $y_{1}w$, $wy_{2}$, $y_{2}x_{n}$, $y_{2}v$, $y_2v_2$,
$y_2p_2$ and $y_2p$ straight. Then the edges $v_2p$, $p_2p$, $v_1v_2$ and $wp_2$ are also straight.

We first color $p$, $v$ and $w$ with $\alpha_{1}$, and set $I_1=\{(p,\alpha_1),(v,\alpha_1),(w,\alpha_1)\}$.
So $|L(v_{2},I_1)|\geq 3$, $|L(y_{2},I_1)|\geq 3$ and $|L(p_{2},I_1)|\geq 3$.
Let $(v,\alpha_{1})(y_{1},\alpha_{2})\in M_{L,vy_{1}}$. If $\alpha_2=\alpha_1$,
then $|L(v,\{y_1,\alpha_1),(y_2,\alpha_1)\}|\ge3$ and $|L(w,\{y_1,\alpha_1),(y_2,\alpha_1)\}|\ge3$.
So there exists the required independent set by Claim~\ref{y1vwy2}. Hence, $\alpha_{1} \neq \alpha_{2}$.
Suppose that $(p,\alpha_{1})(p_{1},\beta) \in M_{L,p_{1}p}$ for some $\beta\in L(p_1)$.
If $\beta \notin\{\alpha_{1},\alpha_{2}\}$, then we first color $y_{1}$ with $\beta$ and set $I_2=I_1\cup\{(y_1,\beta)\}$.
Hence, $|L(x_{n},I_2)|\geq 2$ and $|L(y_{2},I_2)|\geq 3$.
So we can color $y_{2}$ with $\alpha_{3}\in L(y_2,I_2)$ such that
$|L(x_{n},I_2\cup\{(y_2,\alpha_3)\})|\geq 2$.
Note that $|L(p_{1},I_2\cup\{(y_2,\alpha_3)\})|\geq 2$. Now we can color $v_{1}$, $v_{2}$, $z_{1},\ldots, z_{m}$, $p_{1}$, $p_{2}$,
$x_{1},\ldots, x_{n}$, successively.
So suppose that $\beta\in\{\alpha_{1},\alpha_2\}$.
Let $(z_{m},\gamma)(p,\alpha_{1}) \in M_{L,z_{m}p}$.
Then $\gamma\neq\beta$ by the Remark. If $\gamma \notin\{\alpha_1,\alpha_{2}\}$, then we
color $y_{1}$ with $\gamma$ such that $|L(z_{m},I_1\cup\{(y_1,\gamma)\})|\geq 3$.
Hence, $|L(x_{n},I_1\cup\{(y_1,\gamma)\})|\geq 2$ and $|L(y_{2},I_1\cup\{(y_1,\gamma)\})|\geq3$.
So we can color $y_{2}$ with $\mu\in L(y_2,I_1\cup\{(y_1,\gamma)\})$ such that $|L(x_{n},I_1\cup\{(y_1,\gamma),(y_2,\mu)\})|\geq 2$.
Now we can color $p_{1}$, $p_{2}$, $v_{1}$, $v_{2}$, $z_{1},\ldots, z_{m}$, $x_{1},\ldots, x_{n}$, successively.
So suppose that $\gamma\in\{\alpha_1,\alpha_{2}\}\setminus\{\beta\}$.
Recall that $(v_{1},\beta)(p,\alpha_{1}) \notin M_{L,v_{1}p}$ by the Remark.
If $(v_{1},\gamma)(p,\alpha_{1}) \in M_{L,v_{1}p}$, then
$|L(v_{1},\{(p,\alpha_1),(y_1,\gamma)\})|\geq 3$ and $|L(z_{m},\{(p,\alpha_1),(y_1,\gamma)\})|\geq 3$.
By Claim~\ref{y1p1pzm}, there exists the required independent set.
Hence, $(v_{1},\gamma)(p,\alpha_{1}) \notin M_{L,v_{1}p}$, say $(v_{1},\gamma')(p,\alpha_{1}) \in M_{L,v_{1}p}$,
where $\gamma' \notin\{\alpha_{1},\alpha_{2}\}$. We first color $y_{1}$ with $\gamma'$ such that $|L(v_{1},I_1\cup\{(y_1,\gamma')\})|\geq 2$.
Hence, $|L(x_{n},I_1\cup\{(y_1,\gamma')\})|\geq 2$ and $|L(y_{2},I_1\cup\{(y_1,\gamma')\})|\geq 3$.
So we can color $y_{2}$ with $\mu'\in L(y_2,I_1\cup\{(y_1,\gamma')\})$ such that $|L(x_{n},I_1\cup\{(y_1,\gamma'),(y_2,\mu')\})|\geq 2$.
Now we can color $p_{1}$, $p_{2}$, $z_{m},\ldots, z_{1}$, $v_{1}$, $v_{2}$, $x_{1},\ldots, x_{n}$, successively.

\medskip
\noindent {\bf Case 11} $G=G_{9}$.
\medskip

Suppose that $C_{1}=vy_{1}x_{1}v$ does not have property $\mathcal{P}$.
Then we can color $v$ with $\alpha_{1}\in L(v)$ and color $y_{1}$ with $\alpha_{2}\in L(y_{1})$ such that
$|L(x_{1},I')|\geq 3$, where $I'=\{(v,\alpha_1),(y_1,\alpha_2)\}$.
Hence, $|L(v_{1},I')|\geq 2$ and $|L(z_{1},I')|\geq 3$.
So we can color $z_{1}$ with $\alpha_{3}\in  L(z_{1},I')$ such that
$|L(v_{1},I'')|\geq 2$, where $I''=I'\cup\{(z_1,\alpha_3)\}$.
Similarly, we can color $w$ with $\alpha_{4}\in  L(w,I'')$ such that
$|L(p_{1},I''')|\geq 2$, where $I'''=I''\cup\{(w,\alpha_4)\}$.
Now we can color $y_{2}$, $x_{n},\ldots, x_{1}$, $p_{2}$, $p_{1}$, $v_{2}$, $v_{1}$, successively.

So $C_{1}=vy_{1}x_{1}v$ has property $\mathcal{P}$.
By symmetry, $C_{2}=vy_{2}x_{1}v$ has property $\mathcal{P}$, too.
That is, the edges $vy_1$, $y_1x_1$, $x_1v$, $vy_2$ and $y_2x_1$ are straight.
We first color $y_{1}$ and $y_{2}$ with the same color $\alpha_{1}$, and set $I'=\{(y_1,\alpha_1),(y_2,\alpha_1)\}$.
Hence, $|L(w,I')|\geq 2$ and $|L(p_{2},I')|\geq 3$.
So we can color $p_{2}$ with $\alpha_{2}\in  L(p_{2},I')$  such that $|L(w,I'')|\geq 2$, where $I''=I'\cup\{(p_2,\alpha_2)\}$.
Note that $|L(v,I'')|\geq 3$ and $|L(x_{1},I'')|\geq 3$.
Now we can color $z_{1}$, $p_{1}$, $w$, $x_{n},\ldots, x_{2}$, $v_{1}$, $v_{2}$, $v$, $x_{1}$, successively.

\medskip
\noindent{\bf Case 12} $G=G_{10}$.
\medskip

Suppose that $C_{1}=yu_{0}u_{1}y$ does not have property $\mathcal{P}$.
Then we can color $u_{0}$ with $\alpha_{1}\in L(u_{0})$ and color $u_{1}$ with $\alpha_{2}\in L(u_{1})$ such that
$|L(y,I')|\geq 3$, where $I'=\{(u_0,\alpha_1),(u_1,\alpha_2)\}$.
Hence, $|L(u,I')|\geq 2$ and $|L(u_{3},I')|\geq 3$.
So we can color $u_{3}$ with $\alpha_{3}\in  L(u_{3},I')$ such that $|L(u,I'')|\geq 2$, where $I''=I'\cup\{(u_3,\alpha_3)\}$.
Similarly, we can color $z$ with $\alpha_{4}\in  L(z,I'')$ such that
$|L(x,I''')|\geq 2$, where $I'''=I''\cup\{(z,\alpha_4)\}$.
Now we can color $u_{2}$, $u$, $w$, $x$, $y$, successively.

Suppose that $C_{2}=xyu_{0}x$ does not have property $\mathcal{P}$.
Then we can color  $x$ with $\alpha_{1}\in L(x)$
and color $u_{0}$ with $\alpha_{2}\in L(u_0)$ such that $|L(y,I')|\geq 3$,
where $I'=\{(x,\alpha_1),(u_0,\alpha_2)\}$.
Hence, $|L(u_{3},I')|\geq 2$ and $|L(w,I')|\geq 3$.
So we can color $w$ with $\alpha_{3}\in L(w,I')$ such that $|L(u_{3},I'')|\geq 2$, where $I''=I'\cup\{(w,\alpha_3)\}$.
Similarly, we can color $u_{2}$ with $\alpha_{4}\in L(u_2,I'')$ such that
$|L(z,I''')|\geq 2$, where $I'''=I''\cup\{(u_2,\alpha_4)\}$.
Now we can color $u_{3}$, $u$, $u_{1}$, $z$, $y$, successively.

So suppose that both $C_1$ and $C_2$ have property $\mathcal{P}$.
That is, the edges $yu_0$, $u_0u_1$, $u_1y$, $xy$ and $xu_0$ are straight.
We first color $u_{1}$ and $x$ with the same color $\alpha_{1}$, and set $I'=\{(x,\alpha_1),(u_1,\alpha_1)\}$.
Hence, $|L(z,I')|\geq 2$ and $|L(u_{2},I')|\geq 3$.
So we can color $u_{2}$ with $\alpha_{2}\in L(u_2,I')$ such that $|L(z,I'')|\geq 2$,
where $I''=I'\cup\{(u_2,\alpha_2)\}$. Note that $|L(y,I'')|\geq 3$ and $|L(u_{0},I'')|\geq 3$.
Now we can color $u_{3}$, $w$, $z$, $u$, $u_{0}$, $y$, successively.

\begin{figure}
  \centering
  \includegraphics[width=0.90\textwidth]{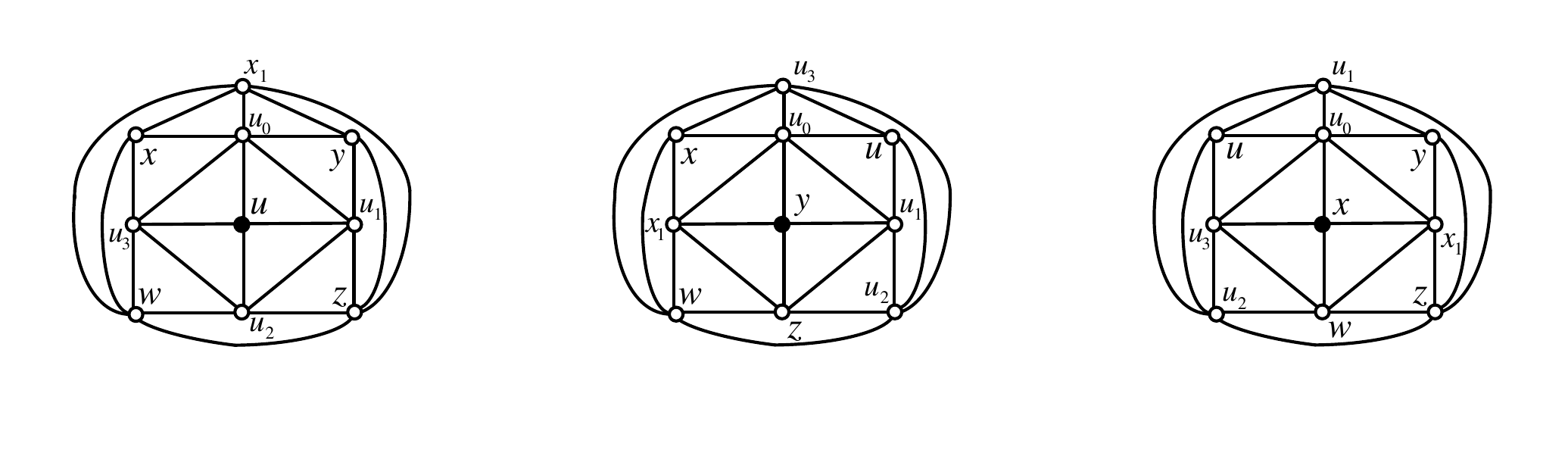}\\
  \vskip-8mm
 $G_{11}$\hspace*{48mm}$G'$\hspace*{48mm}$G''$
  \caption{Two automorphism graphs $G'$ and $G''$ of $G_{11}$.}\label{fig.2}
\end{figure}

\medskip
\noindent{\bf Case 13} $G=G_{11}$.
\medskip

Suppose that $C_{1}=xu_{0}u_{3}x$ does not have property $\mathcal{P}$.
Then we can color $u_{0}$ with $\alpha_{1}\in L(u_0)$ and $u_{3}$ with $\alpha_{2}\in L(u_3)$
such that $|L(x,I')|\geq 3$, where $I'=\{(u_0,\alpha_1),(u_3,\alpha_2)\}$.
Hence, $|L(u,I')|\geq 2$ and $|L(u_{1},I')|\geq 3$.
So we can color $u_{1}$ with $\alpha_{3}\in L(u_1,I')$ such that $|L(u,I'')|\geq 2$, where $I''=I'\cup\{(u_1,\alpha_3)\}$.
Similarly, we can color $z$ with $\alpha_{4}\in L(z,I'')$ such that $|L(y,I''')|\geq 2$, where $I'''=I''\cup\{(z,\alpha_4)\}$.
Now we  color $u_{2}$, $w$, $x_{1}$, $y$, $u$, $x$, successively.

So suppose that $C_{1}=xu_{0}u_{3}x$ has property $\mathcal{P}$.
By the automorphism of $G_{11}$ (as seen in Figure 2), $C_{2}=xu_{0}x_{1}x$ and $C_{3}=uu_{0}u_{3}u$ have property $\mathcal{P}$, too.
By symmetry, $C_{4}=yu_{0}x_{1}y$,
$C_{5}=yu_{0}u_{1}y$ and $C_{6}=uu_{0}u_{1}u$ have property $\mathcal{P}$. That is,
all edges in $C_{1}$, $C_{2},\ldots, C_{6}$ are straight.

\begin{itemize}
\item Suppose that $C_{7}=uu_{3}u_{2}u$ does not have property $\mathcal{P}$.
We can color $u_{3}$ with $\alpha_1$ and color $u_2$ with $\alpha_{2}$ such that
$|L(u,I')|\ge3$, where $I'=\{(u_3,\alpha_1),(u_2,\alpha_2)\}$.
Then we color $x_{1}$ with $\alpha_{1}$, and set $I''=I'\cup\{(x_1,\alpha_1)\}$.
Note that $|L(x,I'')|\geq 3$ and $|L(u_{0},I'')|\geq 3$, since all edges in $C_{1}$, $C_{2},\ldots, C_{6}$ are straight.
Now we color $w$, $z$, $u_{1}$, $y$, $u_{0}$, $u$, $x$, successively.

\item Suppose that  $C_{7}=uu_{3}u_{2}u$ has property $\mathcal{P}$.
By symmetry, $C_{8}=uu_{1}u_{2}u$ has property $\mathcal{P}$, too.
We first color $x_{1}$, $u_{3}$ and $u_{1}$ with the same color $\alpha_{1}$,
and set $I'=\{(x_1,\alpha_1),(u_3,\alpha_1),(u_1,\alpha_1)\}$.
Then $|L(v,I')|\geq 3$ for any vertex $v\in\{x,u_{0},y,u,u_{2}\}$.
Now we color $w$, $z$, $u_{2}$, $y$, $u_{0}$, $u$, $x$, successively.
 \end{itemize}

\medskip
\noindent{\bf Case 14} $G=G_{12}$.
\medskip

Suppose that $C_{1}=yu_{0}u_{1}y$ does not have property $\mathcal{P}$.
Then we color $u_{0}$ with $\alpha_{1}\in L(u_0)$ and $u_{1}$ with $\alpha_{2}\in L(u_1)$
such that $|L(y,I')|\geq 3$, where $I'=\{(u_0,\alpha_1),(u_1,\alpha_2)\}$.
Hence, $|L(u,I')|\geq 2$ and $|L(u_{3},I')|\geq 3$.
So we can color $u_{3}$ with $\alpha_{3}\in L(u_3,I')$ such that $|L(u,I'')|\geq 2$, where $I''=I'\cup\{(u_3,\alpha_3)\}$.
Similarly, we can color $w$ with $\alpha_{4}\in L(w,I'')$ such that $|L(x,I''')|\geq 2$, where $I'''=I''\cup\{(w,\alpha_4)\}$.
Now we color $u_{2}$, $u$, $z$, $y_{1}$, $x_{1}$, $x$, $y$, successively.

Suppose that $C_{2}=uu_{0}u_{1}u$ does not have property $\mathcal{P}$.
Then we can color $u_{0}$ with $\alpha_{1}\in L(u_0)$ and $u_{1}$ with $\alpha_{2}\in L(u_1)$
such that $|L(u,I')|\geq 3$, where $I'=\{(u_0,\alpha_1),(u_1,\alpha_2)\}$.
Hence, $|L(y,I')|\geq 2$ and $|L(y_{1},I')|\geq 3$.
So we can color $y_{1}$ with $\alpha_{3}\in L(y_1,I')$ such that $|L(y,I'')|\geq 2$, where $I''=I'\cup\{(y_1,\alpha_3)\}$.
Similarly, we can color $w$ with $\alpha_{4}\in L(w,I'')$ such that $|L(x_{1},I''')|\geq 2$, where $I'''=I''\cup\{(w,\alpha_4)\}$.
Now we color $z$, $u_{2}$, $u_{3}$, $u$, $x$, $x_{1}$, $y$, successively.

So suppose that both $C_{1}=yu_{0}u_{1}y$ and $C_{2}=uu_{0}u_{1}u$ have property $\mathcal{P}$.
That is, the edges in $C_1$ and $C_2$ are straight. If $C_{3}=yu_{1}y_{1}y$ has property $\mathcal{P}$,
then we can color $u_{0}$ and $y_{1}$ with the same color $\alpha_{1}$, and set $I'=\{(u_0,\alpha_1),(y_1,\alpha_1)\}$.
Hence, $|L(x_{1},I')|\geq 2$ and $|L(w,I')|\geq 3$.
So we can color $w$ with $\alpha_{2}\in L(w,I')$ such that $|L(x_{1},I'')|\geq 2$, where $I''=I'\cup\{(w,\alpha_2)\}$.
Similarly, we can color $u_{2}$ with $\alpha_{3}\in L(u_2,\alpha_3)$ such that
$|L(z,I''')|\geq 2$, where $I'''=I''\cup\{(u_2,\alpha_3)\}$.
Note that $|L(y,I''')|\geq 3$ and $|L(u_{1},I''')|\geq 2$.
Now we can color $u_{3}$, $x$, $x_{1}$, $u$, $u_{1}$, $z$, $y$, successively.
So we may suppose that $C_{3}=yu_{1}y_{1}y$ does not have property $\mathcal{P}$.
That is, we can color $u_1$ with $\alpha_1$ and color $y_1$ with $\alpha_2$ such that
$|L(y,I')|\ge3$, where $I'=\{(u_1,\alpha_1),(y_1,\alpha_2)\}$.
By Lemma~\ref{sraight}, we can rename the colors in $L(u_2)$ to make $uu_2$ straight.

\begin{itemize}
\item If $u_2$ can be colored with $\alpha_1$, then $|L(u,I'')|\geq 3$, where
$I''=I'\cup\{(u_{2},\alpha_1)\}$. Then we color $z$ with $\alpha_3$ and $w$ with $\alpha_4$
 such that $I'''=I''\cup\{(z,\alpha_3),(w,\alpha_4)\}$ is independent set in $H_L$.
 Hence, $|L(u_{3},I''')|\geq 2$ and $|L(x,I''')|\geq 3$.
So we can color $x$ with $\alpha_{5}\in L(x,I''')$ such that $|L(u_{3},I^*)|\geq 2$, where $I^*=I'''\cup\{(x,\alpha_5)\}$.
Now we can color $x_{1}$, $u_{0}$, $y$, $u_{3}$, $u$, successively.

\item Suppose that $u_2$ cannot be colored with $\alpha_1$. That is,
$(u_2,\alpha_1)(u_1,\alpha_1)\in M_{L,u_1u_2}$.
Then we remove the colors of $u_1$ and $y_1$, and color $u_{0}$ and $u_{2}$ with $\alpha_{1}$.
Set $I''=\{(u_0,\alpha_1),(u_2,\alpha_1)\}$.
Note that $|L(u,I'')|\geq 3$ and $|L(u_{1},I'')|\geq 3$.
Since $|L(u_{3},I'')|\geq 2$ and $|L(w,I'')|\geq 3$, we can color $w$ with
$\alpha_{3}\in L(w,I'')$ such that $|L(u_{3},I''')|\geq 2$, where $I'''=I''\cup\{(w,\alpha_3)\}$.
Similarly, we can color $y_{1}$ with $\alpha_{4}\in L(y_1,I''')$ such that $|L(z,I^*)|\geq 2$,
where $I^*=I'''\cup\{(y_1,\alpha_4)\}$. Now we can color $x_{1}$, $y$, $u_{1}$, $z$, $x$,
$u_{3}$, $u$, successively.
\end{itemize}

\medskip
\noindent{\bf Case 15} $G=G_{13}$.
\medskip

Suppose that $C_{1}=zz_{1}u_{2}z$ does not have property $\mathcal{P}$.
Then we can color $z_{1}$ with
$\alpha_{1}\in L(z_1)$ and color $u_{2}$ with $\alpha_{2}\in L(u_2)$ such that $|L(z,I')|\geq 3$, where $I'=\{(z_1,\alpha_1),(u_2,\alpha_2)\}$.
Hence, $|L(w,I')|\geq 2$ and $|L(u_{3},I')|\geq 3$.
So we can color $u_{3}$ with $\alpha_{3}\in L(u_3,I')$ such that $|L(w,I'')|\geq 2$, where $I''=I'\cup\{(u_3,\alpha_3)\}$.
Similarly, we can color $u_{0}$
with $\alpha_{4}\in L(u_0,I'')$ such that $|L(u,I''')|\geq 2$, where $I'''=I''\cup\{(u_0,\alpha_4)\}$.
Now we can color $x$, $w$, $x_{1}$, $y$, $u_{1}$, $u$, $z$, successively.

So suppose that $C_{1}=zz_{1}u_{2}z$ has property $\mathcal{P}$.
By symmetry, $C_{2}=wz_{1}u_{2}w$ has property $\mathcal{P}$, too.
That is, all edges incident with $C_1$ and $C_2$ are straight.
By Lemma~\ref{sraight}, we can rename the colors in $L(u_3)$ and $L(u_1)$
to make the edges $wu_3$ and $zu_1$ straight.

  \begin{itemize}
\item If there exists $\alpha_1$ such that $(u_1,\alpha_1)(u_2,\alpha_1)\in M_{L,u_1u_2}$, then
we color $z_{1}$ and $u_{1}$ with $\alpha_1$, and set $I'=\{(z_1,\alpha_1),(u_1,\alpha_1)\}$.
Then $|L(z,I')|\geq 3$ and  $|L(u_{2},I')|\geq 3$.
Since $|L(y,I')|\geq 2$ and $|L(u_{0},I')|\geq 3$, we can color $u_{0}$ with $\alpha_{2}\in L(u_0,I')$
such that $|L(y,I'')|\geq 2$, where $I''=I'\cup\{(u_0,\alpha_2)\}$.
Similarly, we can color $u_{3}$ with $\alpha_{3}\in L(u_3,I'')$ such that
$|L(u,I''')|\geq 2$, where $I'''=I''\cup\{(u_3,\alpha_3)\}$.
Now we can color $x$, $x_{1}$, $y$, $w$, $u_{2}$, $u$, $z$, successively.

\item So  $(u_1,\alpha)(u_2,\beta)\in M_{L,u_1u_2}$ implies that $\alpha\neq\beta$.
By symmetry, $(u_3,\alpha^*)(u_2,\beta^*)\in M_{L,u_1u_2}$ implies that $\alpha^*\neq\beta^*$.
Now we color $u_{2}$, $u_{1}$ and $u_{3}$ with $\alpha_1$, and set $I'=\{(u_i,\alpha_1) : i=1,2,3\}$.
Then $|L(w,I')|\geq 3$ and $|L(z,I')|\geq 3$.
Next we color $u$ with $\alpha_2\in L(u,I')$ and set $I''=I'\cup\{(u,\alpha_2)\}$;
color  $u_{0}$ with $\alpha_3\in L(u_0,I'')$ and set $I'''=I''\cup\{(u_0,\alpha_3)\}$.
Hence, $|L(x,I''')|\geq 2$, $|L(x_1,I''')|\geq 3$ and $|L(y,I''')|\geq 2$.
So we can color $x$ with $\alpha_{4}\in L(x,I''')$ and color $y$ with $\alpha_{5}\in L(y,I''')$
such that $|L(x_{1},I^*)|\geq 2$, where $I^*=I'''\cup\{(x,\alpha_4),(y,\alpha_5)\}$.
Note that $|L(z,I^*)|\ge2$ and $|L(w,I^*)|\ge2$.
Now we can color $z_{1}$, $x_{1}$, $w$, $z$, successively. \qed
 \end{itemize}

Given a planar graph $G$ of diameter at most two with $|V (G)| \geq 3$, let $G^*$ denote a maximal
planar graph containing $G$ as a spanning subgraph. Obviously, either $G^*=K_{r}$ for $r = 3,4$,
or $G^*$ is a $MP2$-graph. In either case, we have $\chi_{DP}(G^*)\leq 4$ by Theorem~\ref{main2}. By the
fact that $\chi_{DP}(G)\leq \chi_{DP}(G^*)$, we obtain our main theorem.

\begin{theorem}\label{main1}
Every planar graph with diameter at most two is {\rm DP}-$4$-colorable.
\end{theorem}

\end{document}